\documentclass[11pt]{amsart}
\usepackage[margin=1.1in]{geometry}
\usepackage{lineno}
\usepackage{hyperref}

\def\F{\mathcal{F}}
\def\eps{\varepsilon}

\def\d{\delta}

\def\ct{\lceil t \rceil}
\def\clr{\lceil r \rceil}

\newtheorem{theorem}{Theorem}
\newtheorem{corollary}[theorem]{Corollary}
\newtheorem{proposition}[theorem]{Proposition}

\newtheorem{problem}[theorem]{Problem}
\newtheorem{construction}[theorem]{Construction}
\newtheorem{fact}[theorem]{Fact}

\newtheorem{remark}[theorem]{Remark}

\begin{document}

\title{Shadows of 3-uniform hypergraphs under a minimum degree condition}
\thanks{Zolt\'an F\"uredi is partially supported by NKFIH grant KH130371 and NKFI-133819.
Yi Zhao is partially supported by NSF grants DMS 1400073 and 1700622 and Simons Collaboration Grant 710094.}
\author{Zolt\'an F\"uredi}
\address{Alfr\' ed R\' enyi Institute of Mathematics, Budapest, Hungary and Department of Mathematics, University of Illinois at Urbana-Champaign, Urbana, IL, USA} 
\email{z-furedi@illinois.edu} 
\author{Yi Zhao}
\address{Department of Mathematics and Statistics, Georgia State University, Atlanta, GA 30303}
\email{yzhao6@gsu.edu}
\date{\today}

\begin{abstract}
We prove a minimum degree version of the Kruskal--Katona theorem for triple systems: given $d\ge 1/4$ and a triple system $\F$ on $n$ vertices with minimum degree $\d(\F)\ge d\binom n2$, we obtain asymptotically tight lower bounds for the size of its shadow. Equivalently, for $t\ge n/2-1$, we asymptotically determine the minimum size of a graph on $n$ vertices, in which every vertex is contained in at least $\binom t2$ triangles. This can be viewed as  a variant of the Rademacher--Tur\'an problem.
\end{abstract} 

\maketitle

\section{Introduction}
Given a set $X$ and a family $\F$ of $k$-subsets of $X$, the \emph{shadow} $\partial \F$ of $\F$ is the family of all $(k-1)$-subsets of $X$ contained in some member of $\F$. The Kruskal--Katona theorem \cite{Kat68, Kru63} is one of the most important results in extremal set theory -- it gives a tight lower bound for the size of shadows of all $k$-uniform families of a given size. The following is a version due to Lov\'asz \cite{Lov-book}, where $\binom tk = \frac{t(t-1)\cdots (t-k+1)}{k!}$ for a real number $t$. Note that it is tight when $t$ is an integer by considering the family of all $k$-subsets of a set of $t$ vertices.

\begin{theorem}[Kruskal--Katona theorem]
\label{thm:KK}
If $\F$ is a family of $k$-sets with $|\F|\ge \binom{t}{k}$ for some real number $t$, then $|\partial \F|\ge \binom{t}{k-1}$.
\end{theorem}

A family $\F$ of $k$-subsets of $X$ is often regarded as a \emph{$k$-uniform hypergraph}, or \emph{$k$-graph} $(X, \F)$ with $X$ as the vertex set and $\F$ as the edge set. 
For every $x\in X$, define $\F_x = \{F\setminus x: x\in F$ and $F\in \F\}$. The minimum (vertex) degree of $\F$ is denoted by $\d(\F):= \min_{x} |\F_x|$. 
The following minimum degree version of the Kruskal--Katona theorem has not been studied before but emerged naturally when Han, Zang, and Zhao \cite{HZZ17} investigated a packing problem for 3-graphs.

\begin{problem}
\label{pro:KK}
Given $k\ge 3$ and $0<d<1$, let $X$ be a set of $n$ vertices and $\F$ be a family of $k$-subsets of $X$ with $\d(\F)\ge d \binom{n}{k-1}$.\footnote{ It is more natural to assume $\d(\F)\ge d \binom{n-1}{k-1}$ as $\binom{n-1}{k-1}$ is the largest possible degree. However, since we are mainly interested in the asymptotics of $|\partial \F|$, we choose the simpler looking condition $\d(\F)\ge d \binom{n}{k-1}$.}
How small can $|\partial \F|$ be?
\end{problem}

Problem~\ref{pro:KK} belongs to an area of active research on extremal problems under maximum or minimum degree conditions.
Two early examples are the work of Bollob{\'a}s, Daykin, and Erd{\H{o}}s \cite{BDE}, who studied the minimum degree version of the Erd\H{o}s matching conjecture, and of Frankl \cite{Fra87}, who studied the Erd\H{o}s--Ko--Rado theorem under maximum degree conditions. More recent examples include the minimum (co)degree Tur\'an's problems \cite{LoMa-de, MZ07b}, the minimum degree version of the Erd\H{o}s--Ko--Rado theorem \cite{FrTo16, HuZh16, MR3974705}, and the minimum degree version of Hilton--Milner theorem \cite{MR3741439, MR3974705}. Recently Jung \cite{Jung} studied minimum $|\partial \F|/ |\F|$ among all $k$-graphs $\F$ with maximum degree $\Delta(\F)\le d$.

\smallskip

Since $\d(\F)\ge d \binom{n}{k-1}$ implies that $|\F|\ge d\binom nk$, we could apply Theorem~\ref{thm:KK} to $\F$ but will not obtain a tight bound for $|\partial \F|$. A better approach is applying Theorem~\ref{thm:KK}  to $\F_x$ for each vertex $x$. Since $|\F_x| \ge d \binom{n}{k-1} \ge \binom{d^{\frac1{k-1}}n}{k-1}$, by Theorem~\ref{thm:KK}, we have $|\partial \F_x| \ge \binom{d^{\frac1{k-1}}n}{k-2} \ge d^{\frac{k-2}{k-1}} \binom{n}{k-2} + O(n^{k-3})$. Consequently,
\begin{align}\label{eq:dF2}
|\partial \F| = \sum_{x} \frac{|\partial \F_x|}{k-1} \ge \frac{n}{k-1} d^{\frac{k-2}{k-1}} \binom{n}{k-2} + O(n^{k-2})
\ge d^{\frac{k-2}{k-1}} \binom{n}{k-1} + O(n^{k-2}).
\end{align}
This bound is tight (up to the error term) when the first inequality in \eqref{eq:dF2} is asymptotically an equality, which occurs when $\F_x$ is a clique of order $d^{\frac1{k-1}}n$ for every $x$. Thus, the bound in \eqref{eq:dF2} is asymptotically tight when $\F$ consists of $d^{\frac{1}{1-k}}$ vertex-disjoint cliques of order $d^{\frac1{k-1}}n$, in particular, when $d= \ell^{1-k}$ for some $\ell\in \mathbb{N}$. 

In this paper we improve \eqref{eq:dF2} and answer Problem~\ref{pro:KK} asymptotically for $k=3$ and $d\ge 1/4$. Two overlapping cliques of order about $\sqrt{d}n+1$ is a natural candidate for extremal hypergraphs -- the following theorem confirms this for $\frac14 \le d< \frac{47 - 5\sqrt{57}}{24} \approx 0.385$. However, there is a different extremal hypergraph for larger values of $d$.
\begin{theorem}
\label{thm:main}
Let $1/4 \le d <1$ and $n\in \mathbb{N}$ be sufficiently large. If $\F$ is a triple system on $n$ vertices with $\d(\F)\ge d \binom{n}2$, 
then
\[
|\partial \F| \ge \left\{
\begin{array}{ll}
\left(4\sqrt{d}-2d-1\right) \binom n2 & \text{if  } \frac14 \le d< \frac{47 - 5\sqrt{57}}{24} \\
 \left( \frac12+ \sqrt{\frac{4d - 1}{12}} \right)\binom n2 & \text{if  } d\ge \frac{47 - 5\sqrt{57}}{24} .
\end{array} \right.
\]
These bounds are best possible up to an additive term of $O(n)$.
\end{theorem}

Although seemingly technical, Theorem~\ref{thm:main} has an interesting application on \emph{3-graph packing and covering}. Given positive integers $a, b, c$, let $K^3_{a, b, c}$ denote the complete 3-partite \emph{3-graph} with parts of size $a, b$, and $c$. Answering a question of Mycroft \cite{MR3423473}, Han, Zang, and Zhao \cite{HZZ17} determined the minimum $\d(H)$ of a 3-graph $H$ that forces a \emph{perfect $K^3_{a, b, c}$-packing} in $H$ for any given $a, b, c$.\footnote{ Given hypergraphs $H$ and $F$, a perfect $F$-packing in $H$ is a spanning subgraph of $H$ that consists of vertex-disjoint copies of $F$.} One of the main steps in their proof is determining the smallest $\d(H)$ of a 3-graph $H$ that guarantees that every vertex of $H$ is contained in a copy of $K^3_{a, b, c}$ (this is necessary for $H$ containing a perfect $K^3_{a, b, c}$-packing). 

\begin{corollary} \cite[Lemma 3.7]{HZZ17}
\label{cor:HZZ}
Let $d_0= 6- 4\sqrt{2} \approx 0.343$.
For any $\gamma > 0$, there exists $\eta > 0$ such that the following holds for sufficiently large $n$. If $H$ is an $n$-vertex 3-graph with $\delta(H) \ge (d_0 + \gamma) \binom n2$, then each vertex of $H$ is contained in at least $\eta n^{a+b+c-1}$ copies of $K^3_{a, b, c}$.
\end{corollary}

It was shown \cite[Construction 2.6]{HZZ17} that $d_0$ in Corollary~\ref{cor:HZZ} is best possible. 
We give a proof outline of Corollary~\ref{cor:HZZ} at the end of Section 2 -- a complete proof can be found in \cite{HZZ17}.

\bigskip
Our approach towards Theorem~\ref{thm:main} is viewing it as an extremal problem on graphs. The following is an equivalent form of Problem~\ref{pro:KK}, in which $K_t^k$ denotes the complete $k$-graph on $t$ vertices (and we omit the superscript when $k=2$).
\begin{problem}\label{pro:KK2}
Given a $(k-1)$-graph $G$ on $n$ vertices such that every vertex is contained in at least $d \binom{n}{k-1}$ copies of $K_k^{k-1}$, how many edges must $G$ have?
\end{problem}
To see why Problems~\ref{pro:KK} and \ref{pro:KK2} are equivalent, let $m_1$ be the minimum $|\partial \F|$ for Problem~\ref{pro:KK} and $m_2$ be the minimum $e(G)$ for Problem~\ref{pro:KK2}. To see why $m_1\ge m_2$, consider a $k$-uniform family $\F$ with $\d(\F)\ge d \binom{n}{k-1}$. Let $G= (V(\F), \partial \F)$ be the $(k-1)$-graph of its shadow. Since every member of $\F$ gives rise to a copy of $K_k^{k-1}$ in $G$, $\d(\F)\ge d \binom{n}{k-1}$ implies that every vertex is contained in at least $d \binom{n}{k-1}$ copies of $K_k^{k-1}$. Thus $|\partial \F| = e(G) \ge m_2$.  To see why $m_2\ge m_1$, consider a $(k-1)$-graph $G$ such that every vertex is contained in at least $d \binom{n}{k-1}$ copies of $K_k^{k-1}$. Let $\F$ be the family of $k$-subsets of $V(G)$ that span copies of $K_k^{k-1}$ in $G$. Then $\partial F \subseteq G$ and for every $v\in V(G)$, we have $|\F_v| \ge d \binom{n}{k-1}$. Thus $e(G)\ge |\partial \F|\ge m_1$ as desired.

\medskip
In order to prove Theorem~\ref{thm:main}, we solve the $k=3$ case of Problem~\ref{pro:KK2} with $d\ge 1/4$. 
For convenience, we assume that every vertex of $G$ is contained in at least $\binom t2$ triangles. 
There are essentially two \emph{extremal graphs}: the first one consists of two copies of $K_{t+1}$ that share $2t+2 - n$ vertices; the second one is obtained from two disjoint copies of $K_{n/2}$ by adding a regular bipartite graph between them. The size of these two extremal graphs can be conveniently represented by a quadratic function $f(x)$, which arises naturally from a lower bound for $e(G)$ in Proposition~\ref{clm:eG}. 
\begin{theorem} \label{thm:G}
Let $n\in \mathbb{N}$, $t, r\in \mathbb{R}$ such that $n/2 \le t+1\le n$, $r\ge 0$, and 
\begin{equation}
\label{eq:r}
\binom{\frac{n}{2}-1}{2} + 3\binom{r}{2} = \binom{t}{2}.
\end{equation}
Define a function $f: \mathbb{R} \to \mathbb{R}$ as  
\begin{align}\label{eq:f}
f(x) = \binom{t}{2} + x (n-x) - \binom{n-x -1}{2}. 
\end{align}
If $G$ is an $n$-vertex graph such that each vertex is contained in at least $\binom t2$ triangles, then
\begin{align}\label{eq:eG}
e(G) \ge \left\{
\begin{array}{ll}
f(t) & \text{if  } r+t\le \frac{5n}6 \text{ or approximately } t\le 0.6208 n,\\ 
f(\frac{n}{2}+r-1) & \text{otherwise.}
\end{array} \right.
\end{align}
Furthermore, these bounds are tight when $n/2, t, r$ are integers, and tight up to an additive $O(n)$ in general. 
\end{theorem}

Theorem~\ref{thm:G} can be viewed as a variant of the well-studied Rademacher--Tur\'an problem. Starting with the work of Rademacher (unpublished) and of Erd\H{o}s \cite{Erd55}, the Rademacher--Tur\'an problem studies the minimum number of triangles in a graph with given order and size. Instead of the total number of triangles in a graph, one may ask for the maximum or minimum number of triangles containing a fixed vertex. Given a graph $G$, we define the \emph{triangle-degree} of a vertex as the number of triangles that contain this vertex. Let $\Delta_{K_3}(G)$ and $\delta_{K_3}(G)$ denote the maximum and minimum triangle-degree in $G$, respectively. The contrapositive of Theorem~\ref{thm:G} states that \emph{if $G$ is a graph on $n$ vertices that fails \eqref{eq:eG}, then $\delta_{K_3}(G)< \binom t2$.} Correspondingly, the \emph{maximum} triangle-degree version of Rademacher--Tur\'an problem was recently studied by Falgas-Ravry, Markstr\"om, and Zhao \cite{MR4225783}. In addition, Theorem~\ref{thm:G} looks similar to the question of Erd\H{o}s and Rothschild \cite{MR891250} on the \emph{book size} of graphs: in the complementary form, it asks for the maximum size of a graph on $n$ vertices, in which every edge is contained in at most $d$ triangles.

\medskip

We prove Theorem~\ref{thm:G} and Theorem~\ref{thm:main} in the next section.
When $t < n/2-1$, it is reasonable to speculate that an extremal graph is a disjoint union of copies of $K_{t+1}$ and an extremal graph for Theorem~\ref{thm:G}. Unfortunately we cannot verify this. We provide some evidence for this speculation in the last section. 


\medskip

\noindent \textbf{Notation.}
Given a family $\F$ of sets, $|\F|$ is the size of $\F$, namely, the number of sets in $\F$. 
A $k$-uniform hypergraph $H$, or $k$-graph, consists of a vertex set $V(H)$ and an edge set $E(H)$, which is a family of $k$-subsets of $V(H)$.  
Given a vertex set $S$, denote by $e_H(S)$ the number of edges of $H$ induced on $S$.
Suppose $G$ is a graph. 
For a vertex $v\in V(G)$, let $N_G(v)$ denote the \emph{neighborhood} of $v$, the set of vertices adjacent to $v$, and let $d_G(v)=| N_G(v)|$ be the degree of $v$. Let $N_G[v]:= N_G(v)\cup \{v\}$ denote the \emph{closed neighborhood} of $v$. When the underlying (hyper)graph is clear from the context, we omit the subscript in these notations. 


 \section{Proofs of Theorem~\ref{thm:G} and Thereom~\ref{thm:main} }
Suppose that $G=(V, E)$ is a graph on $n$ vertices such that each vertex is contained in at least $\binom t2$ triangles, in other words,
\begin{equation}\label{eq:N}
\forall v\in V, \quad e(N(v))\ge \binom{t}{2}, 
\end{equation}
where $t$ is a positive real number. 
Trivially $t\le \d(G)\le n-1$ because $e(N(v))\le \binom{d(v)}{2}$ for every vertex $v\in V$. Therefore 
\[
e(G)\ge \frac{\d(G)  n}2 \ge \frac{t n}{2}. 
\]
When $t+1$ divides $n$, this bound is tight because $G$ can be a disjoint union of $\frac{n}{t+1}$ copies of $K_{t+1}$. Below we often assume that $t\le n-2$ because when $t= n-1$, we must have $G=K_n$. 

Let us derive another lower bound for $e(G)$ by using the function $f$ defined in \eqref{eq:f}.

\begin{proposition}\label{clm:eG}
If $G=(V, E)$ is a graph on $n$ vertices satisfying \eqref{eq:N}, then $e(G)\ge f(\d(G))$, and the equality holds if and only if there exists $v_0\in V$ such that $e(N(v_0))= \binom t2$,  $d(v)= \d(G)$ for all $v\not\in N(v_0)$, and $V\setminus N[v_0]$ induces a clique.
\end{proposition}

\begin{proof}
Suppose  $\d(G) = \d$ and $v_0\in V$ satisfies $d(v_0)= \d$. Since we may partition $E(G)$ into the edges induced on $N(v_0)$ and the edges incident to some vertex $v\not\in N(v_0)$, we have  
\[
e(G) = e(N(v_0)) + \left(\sum_{v\not\in N(v_0)} d(v)\right) - e(V\setminus N(v_0)). 
\]
Because of \eqref{eq:N}, $d(v)\ge \delta$ for all $v\not\in N(v_0)$, and $e(V\setminus N(v_0))\le \binom{n-\d -1}{2}$ (note that  $v_0$ has no neighbor outside $N(v_0)$), we derive that $e(G) \ge  \binom{t}{2} + \d (n-\d) - \binom{n-\d -1}{2}$. Furthermore, the equality holds exactly when $e(N(v_0))= \binom t2$,  $d(v)= \d(G)$ for all $v\not\in N(v_0)$, and $V\setminus N[v_0]$ induces a clique.
\end{proof}

Let us construct three graphs satisfying \eqref{eq:N}.  
Note that, if $r$ satisfies \eqref{eq:r}, then $r\le n/2$ because $\binom{ {n}/{2}-1}{2} + 3\binom{{n}/2}{2} = \binom{n-1}{2} \ge \binom t2$.

\begin{construction} \label{cons}
Suppose $t, r\in \mathbb{R}$ satisfy $\frac n2 -1 \le t\le n-2$, $r\ge 0$, and \eqref{eq:r}.
\begin{enumerate}
\item Let $G_1$ be the union of two copies of $K_{\lceil t\rceil +1}$ sharing $2\lceil t \rceil+2 - n$ vertices. 
\item When $n$ is even, let $G_2$ be the $n$-vertex graph obtained from two disjoint copies of $K_{n/2}$ by adding an $\lceil r \rceil$-regular bipartite graph between two cliques. 
\item When $n$ is odd, let $r' \in \mathbb{R}^+$ satisfy $\binom{\frac{n-3}{2}}{2} + 3\binom{r'}{2} = \binom{t}{2}$. 
Let $G'_2$ be the $n$-vertex graph obtained from two disjoint copies of $K_{(n-1)/2}$ by adding an $\lceil r' \rceil$-regular bipartite graph between them, and a new vertex whose adjacency is the exactly the same as one of the existing vertices. 
\end{enumerate}

\end{construction}
It is easy to see that $G_1, G_2, G'_2$ all satisfy \eqref{eq:N}.  For example, consider a vertex $x\in V(G_2)$. Let $A$ and $B$ denote the vertex sets of the two copies of $K_{n/2}$ of $G_2$ and assume $x\in A$. Then $N(x)$ contains $\binom{n/2 -1}2$ edges from $A$, $\binom{\clr}2$ edges from $B$, and $\clr (\clr -1)$ edges between $A$ and $B$. Hence $e(N(x))= \binom{\frac{n}{2}-1}{2} + 3\binom{\clr}{2} \ge \binom{t}{2}$.

\smallskip
The following proposition gives the sizes of $G_1, G_2$, and $G'_2$. 
\begin{proposition}\label{pro:sizeG}
Suppose $n\in \mathbb{N}$, $t, r\ge 0$ satisfy $\frac n2 -1 \le t\le n-1$ and \eqref{eq:r}. 
If all $n/2, t, r$ are integers, then $e(G_1)= f(t)$ and $e(G_2)= f(n/2 +r -1)$, otherwise $e(G_1)\le f(t) + n$ and $e(G_2)\le f(n/2 + r -1) + n/2$. Furthermore, $e(G'_2)=  f(n/2 + r -1) + O(n)$ when $r', r = \Omega(n)$.
\end{proposition}
\begin{proof}
First, by the definition of $f(x)$, it is easy to see that 
\begin{align} \label{eq:fk}
f(t) = \binom n2 - (n-1-t)^2 
\end{align}
(alternatively when $t\in \mathbb{Z}$, we can apply Proposition~\ref{clm:eG} by letting $v_0$ be any vertex not in the intersection of the two cliques). 
We know that 
\[
e(G_1)= \binom n2 - (n-1- \lceil t \rceil)^2 \ge  \binom n2 - (n-1-t)^2 = f(t)
\]
and equality holds when $t\in \mathbb{Z}$. In addition, we have $e(G_1)\le f(t) + n$ because
\begin{align*} 
 (n-1- \lceil t \rceil)^2 - (n-1-t)^2 &= \left( 2(n-1) - (\ct + t ) \right) (\ct - t) \le n
\end{align*}
by using $t+1\ge \ct\ge t\ge n/2 - 1$.

Second, using the definitions of $f(x)$ and $r$, it is not hard to see that 
\begin{align} \label{eq:G2}
f\left(\frac{n}{2}+r-1\right) = \frac{n}{2} \left(\frac{n}{2}+ r -1 \right). 
\end{align}
It follows that 
\[
e(G_2) = \frac{n}{2} \left(\frac{n}{2}+ \lceil r \rceil -1 \right)\le f\left(\frac{n}{2}+r-1\right) + \frac n2
\]
and equality holds when $r\in \mathbb{Z}$.

Third, it is easy to see that 
\[
e(G'_2) = \frac{n+1}2 \left( \frac{n-1}2 + \lceil r' \rceil -1 \right).
\]
By the definitions of $r$ and $r'$, we have $\binom {r'}2 - \binom{r}2 = \frac{2n - 7}{24}$. When $r, r'= \Omega(n)$, we have $r' - r = O(1)$ and consequently,
\begin{align*}
e(G'_2) - f\left(\frac{n}{2}+r-1\right) 
&\le \frac{n+1}2 \left( \frac{n-1}2 + r' -1 \right) - \frac{n}{2} \left(\frac{n}{2}+ r -1 \right) \\
& = \frac n2 (r' -r) + \frac{r'}2 - \frac 34 = O(n). \qedhere
\end{align*}

\end{proof}

We compare $f(t)$, the approximate size of $G_1$, with $f(\frac{n}{2}+r-1)$, the approximate size of $G_2$ and $G'_2$, in the next proposition.
\begin{proposition}\label{pro:com}
Suppose $\frac n2 -1 \le t\le n-1$, $f(x)$ and $r$ are defined as in \eqref{eq:f} and \eqref{eq:r}, respectively. We have $f(t)\le f(\frac{n}{2}+r-1)$ if and only if $r+t\le \frac{5n}6$, equivalently,
\begin{align} \label{eq:t}
 t\le \frac54 n - \frac{\sqrt{ 57 n^2 - 72n}}{12} - 1 \approx 0.6208 n.
\end{align}
\end{proposition}

To prove Proposition~\ref{pro:com}, we need a simple fact on quadratic functions.
\begin{fact} \label{fact}
Suppose $g(x)$ is a quadratic function with a maximum at $x= a$ and assume $x_1 \le x_2$. Then $g(x_1)\le g(x_2)$ if and only if $x_1 + x_2 \le 2a$. \qed
\end{fact}

\begin{proof}[Proof of Proposition~\ref{pro:com}]
First note that 
\[
f(x) = - \frac 32 x^2 + \frac{4n-3}{2} x - \frac{n^2}2 + \binom{t}{2} + \frac32 n - 1
\]
is a quadratic function with a maximum at $x= \frac{2n}{3}- \frac{1}{2}$. Second, since $r\le \frac{n}2$, it follows that
\[
\binom{\frac{n}{2}+r-1}{2} = \binom{\frac n2 -1}2 + \left(\frac n2 - 1 \right) r + \binom r2 \ge \binom{\frac{n}{2}-1}{2} + 3\binom{r}{2} = \binom t2.
\]  
Consequently $\frac{n}{2}+r-1\ge t$. By Fact~\ref{fact}, $f(t)\le f(\frac{n}{2}+r-1)$ if and only if $t+ \frac{n}{2}+r-1 \le \frac{4n}{3}-1$ or $r+ t\le \frac{5n}6$. By \eqref{eq:r}, this is equivalent to 
\[
\binom{\frac{n}{2}-1}{2} + 3\binom{\frac{5n}6 - t}{2} \ge \binom{t}{2} \quad \text{or} \quad 
 (t+1)^2 - \frac{5}{2}(t+1)n + \frac{7}{6} n^2 + \frac{n}{2} \ge 0,
\]
which holds exactly when  $t+1\le \frac54 n - \frac{\sqrt{ 57 n^2 - 72n}}{12}$ (because $t< n$).
\end{proof}

We are ready to prove Theorem~\ref{thm:G}.
\begin{proof}[Proof of Theorem~\ref{thm:G}]
Assume that $\d= \d(G)$. We separate two cases.

\medskip
\noindent \textbf{Case 1:} $r+t\le \frac{5n}6 $, equivalently, \eqref{eq:t}. 

First assume that $\d\ge \frac43 n - t -1$. Since $t\le \frac{5n}6 - r$, we have $\d\ge \frac n2 + r -1$ and consequently,
\[
e(G) \ge \frac{n}2 \left(\frac n2 + r -1 \right) = f\left(\frac n2 + r -1\right)\ge f(t)
\]
by \eqref{eq:G2} and Proposition~\ref{pro:com}.

Second assume that $\d< \frac43 n - t -1$. By Proposition~\ref{clm:eG}, we have $e(G)\ge f(\d)$. Recall that \eqref{eq:N} forces $t\le \d$. Since $t\le \d< \frac43 n - t -1$ and $f(x)$ is a quadratic function maximized at $\frac{2n}{3}- \frac{1}{2}$, we derive from Fact~\ref{fact} that $f(\d)\ge f(t)$. 
Hence $e(G)\ge f(\d) \ge f(t)$.

\medskip
\noindent \textbf{Case 2:} $r+t > \frac{5n}6$. 

If $\d \ge \frac{n}{2}+r-1$, then $e(G)\ge \frac{n}{2} (\frac{n}{2}+r-1)= f(\frac{n}{2}+r-1)$ by \eqref{eq:G2}. Otherwise $\d< \frac{n}{2}+r-1$.
Note that 
\[
\d +  \frac{n}{2}+r-1 \ge t +  \frac{n}{2}+r-1 > \frac{5n}{6} + \frac{n}{2}-1 = \frac{4n}{3}-1.
\]
Since the quadratic function $f(x)$ is maximized at $\frac{2n}{3}- \frac{1}{2}$, we derive from Fact~\ref{fact} that $f(\d) \ge f(\frac{n}{2}+r-1)$. By Proposition~\ref{clm:eG}, we have $e(G)\ge f(\d)\ge f(\frac{n}{2}+r-1)$.

\medskip
By Proposition~\ref{pro:sizeG}, when $n/2, t, r$ are all integers, we have $e(G_1)= f(t)$ and $e(G_2)= f(\frac{n}{2}+r-1)$. In other cases, we have $e(G_1)\le f(t) + n$ and $e(G_2)\le f(\frac{n}{2}+r-1)+ n/2$. When $n$ is odd and $r+t> 5n/6$, we have $r, r'= \Omega(n)$ and thus $e(G'_2)= f(\frac{n}{2}+r-1)+ O(n)$.
\end{proof}

\begin{remark}
When $n/2, t, r$ are all integers, we actually learn the following about extremal graphs from the proof of Theorem~\ref{thm:G}. Suppose that $G$ is an extremal graph. We claim that $G= G_1$ when $r+t< 5n/6$, and $G$ is $(n/2 + r -1)$-regular when $r+t > 5n/6$,.

Indeed, first assume $r+t< 5n/6$. If $\d\ge \frac43 n - t -1$, then $\d > \frac n2 + r -1$ and consequently, $e(G) > \frac{n}{2} (\frac{n}{2}+r-1) = f(t)$, a contradiction. Following the second case of Case~1, we obtain that $e(G)=f(\delta)=f(t)$ and consequently, $\d=t$. Using Proposition~\ref{clm:eG}, we can derive that $G= G_1$. 
When $r+t > 5n/6$, the second case of Case~2 shows that $e(G)\ge f(\d) > f(\frac{n}{2}+r-1)$, a contradiction. Thus $\d \ge \frac{n}{2}+r-1$ and $e(G)= \frac{n}{2} (\frac{n}{2}+r-1)$, which forces $G$ to be $(n/2 + r -1)$-regular. 
\end{remark}

We now prove Theorem~\ref{thm:main} by applying Theorem~\ref{thm:G} and the arguments that show the equivalence of Problems~\ref{pro:KK} and \ref{pro:KK2} in Section 1.


%
%
%


\begin{proof}[Proof of Theorem~\ref{thm:main}]
Suppose $1/4 \le d <1$ and $n\in \mathbb{N}$ is sufficiently large. 
Choose $t\in \mathbb{R}^+$ such that $\binom t2 = d\binom n2$. Since $\binom{\sqrt{d} n}2 < d\binom n2 < \binom{\sqrt{d} n + 1}2$, we have $\sqrt{d} n < t <\sqrt{d} n + 1$.

Suppose $\F$ is a triple system on $n$ vertices with $\d(\F)\ge d \binom{n}2$.
Let $G = (V(\F), \partial \F)$ be the graph whose edge set is the shadow $\partial \F$. For every $x\in V(G)$, we have $e_G( N(x) )\ge d \binom{n}{2}$. 

\smallskip
\noindent \textbf{Case 1:} $\frac14\le d< \frac{47 - 5\sqrt{57}}{24}$.

Thus $\frac12 \le \sqrt{d} < \frac{15- \sqrt{57}}{12}$. Since $n$ is sufficiently large, we have $\sqrt{d}n \le  \frac{15- \sqrt{57}}{12}n - 2$. Since $\sqrt{d} n < t < \sqrt{d} n + 1$, it follows that
\[
\frac{n}2 < t < \frac{15- \sqrt{57}}{12} n - 1 <  \frac54 n - \frac{\sqrt{ 57 n^2 - 72n}}{12} - 1.
\]
This allows us to apply the first case of Theorem~\ref{thm:G} and \eqref{eq:fk} to derive that
\begin{align*}
e(G) &\ge  f(t) = \binom n2 - (n-1-t)^2 \ge  \binom n2 - (n-1- \sqrt{d} n)^2 \\
&= (4\sqrt{d}-2d-1)\binom n2 + n - dn - 1 \\
&\ge (4\sqrt{d}-2d-1)\binom n2 \quad \text{as $d<1$ and $n$ is sufficiently large.}
\end{align*}

\noindent \textbf{Case 2:} $d\ge \frac{47 - 5\sqrt{57}}{24}$.

Thus $\sqrt{d} \ge \frac{15- \sqrt{57}}{12}$. Since $t > \sqrt{d} n$, it follows that 
\[
t+1 > \frac{15- \sqrt{57}}{12} n+1 > \frac54 n -  \frac{\sqrt{57n^2 - 72}}{12} 
\]
because $\sqrt{57n^2 - 72} > \sqrt{57 n^2} - 6$ for $n\ge 2$. Since \eqref{eq:t} fails, we will apply the second case of Theorem~\ref{thm:G}. Since $\binom t2 = d\binom n2$ and $r\ge 0$, we can obtain from \eqref{eq:r} that 
\[
r= \frac16 \left(3+ \sqrt{3(n-1) \big( (4d-1)n + 5\big) }\right) =  \frac 12 + \frac{n}{2} \sqrt{ \frac{4d-1}{3} } + h(n),
\]
where 
\[ h(n) = \frac1{2\sqrt{3}} \left( \sqrt{ (4d-1)n^2 + (6- 4d)n -5 } - \sqrt{4d-1} n \right).
\]
It is easy to see that $0\le h(n)= O(1)$.
Theorem~\ref{thm:G} thus gives that
\begin{align}
e(G) &\ge f\left( \frac n2 + r -1 \right) = \frac n2 \left( \frac n2 + r -1 \right) \nonumber \\
& = \frac n2 \left( \frac n2 -\frac12  + \frac{n}2\sqrt{\frac{4d-1}3} + h(n) \right)  \nonumber \\
&= \binom n2 \left( \frac12+ \sqrt{\frac{4d-1}{12}} \right) + \frac n4 \sqrt{ \frac{4d-1}3 } + \frac n2 h(n) \label{eq:fr}\\
& \ge \binom n2 \left( \frac12+ \sqrt{\frac{4d-1}{12}} \right). \nonumber
\end{align}

\medskip
To see why these bounds are asymptotically tight, for every graph $G\in \{G_1, G_2, G'_2\}$, we construct a triple system $\F_G$ whose members are all triangles of $G$. Then $\partial \F_G \subseteq E(G)$ and $\d(\F_G)\ge \binom t2= d\binom n2$. 

Proposition~\ref{pro:sizeG} gives that $| \partial \F_{G_1}| \le e(G_1)\le f(t) + n$. By \eqref{eq:fk} and the assumption $t \le \sqrt{d} n + 1$, 
\begin{align*}
| \partial \F_{G_1}| & \le f(t) + n \le  \binom n2 - (n-2- \sqrt{d} n)^2 + n  \\
&=  \left(4\sqrt{d}-2d-1 \right)\binom n2 + \left(3- 2\sqrt{d} - d \right) n - 4 + n\\
&= \left(4\sqrt{d}-2d-1 \right)\binom n2 + O(n).
\end{align*}
When $n$ is even, we apply Proposition~\ref{pro:sizeG} and \eqref{eq:fr} obtaining that
\[
| \partial \F_{G_2}| \le e(G_2) \le  f\left( \frac n2 + r -1 \right) + \frac n2 = \binom n2 \left( \frac12+ \sqrt{\frac{4d-1}{12}} \right) + O(n).
\]
When $n$ is odd, we assume $r+t> 5n/6$ and thus $r, r'= \Omega(n)$. By Proposition~\ref{pro:sizeG} and \eqref{eq:fr}, we conclude that
\[
| \partial \F_{G'_2}| \le e(G'_2) = f\left(\frac{n}{2}+r-1\right)+ O(n) = \binom n2 \left( \frac12+ \sqrt{\frac{4d-1}{12}} \right) + O(n). \qedhere
\]

\end{proof}

We outline the proof of Corollary~\ref{cor:HZZ} emphasizing how Theorem~\ref{thm:main} is applied. In a 3-graph, the \emph{degree} of a pair $p$ of vertices is the number of the edges that contains $p$.

\begin{proof}[Proof Outline of Corollary~\ref{cor:HZZ}]

Assume $\eta \ll \gamma$ and $\eps = \gamma/12$. 
Let $H$ be an $n$-vertex 3-graph and $x$ be a vertex of $H$. 
In order to find $\eta n^{a+b+c-1}$ copies of $K_{a, b, c}^3$, it suffices to find $\frac{\gamma}{2} \binom n2$ pairs of vertices of $H_x$ with degree at lease $\eps^2 n$ -- this follows from standard counting arguments in extremal (hyper)graph theory, or conveniently \cite[Lemma 4.2]{LoMa-fa} of Lo and Markstr\"{o}m. 

Suppose $\d_1(H)\ge (d_0 + \gamma) \binom n2$ with  $d_0= 6- 4\sqrt{2} \approx 0.343$. As shown in \cite[Lemma 3.3]{HZZ17}, it is easy to find a set $V_0$ of at most $3\eps n$ vertices and a subgraph $H'$ of $H$ on $V\setminus V_0$ such that $\delta(H') \ge d_0 \binom{n'}2$, where $n' = |V\setminus V_0|$, and every pair in $\partial H'$ has degree at least $\eps^2 n$ in $H$. Since $\frac14 < d_0 < \frac{47 - 5\sqrt{57}}{24} \approx 0.385$, by the first case of Theorem~\ref{thm:main}, we have 
\[
|\partial H'|\ge (4\sqrt{d_0}-2d_0-1) \binom{n'}2 \ge \left(4\sqrt{d_0}-2d_0-1 -\frac{\gamma}2 \right) \binom{n}2.
\] 
For every vertex $x\in V(H)$, since $d(x)\ge (d_0 + \gamma) \binom n2$ and crucially $4\sqrt{d_0}-2d_0-1 = 1 - d_0$, at least $\frac{\gamma}{2} \binom n2$ pairs in $H_x$ are also in $\partial H'$ thus having degree at lease $\eps^2 n$, as desired.

\end{proof}

\section{Concluding remarks}
Let us restate the $k=3$ case of Problem~\ref{pro:KK2}.
\begin{problem}\label{pro:KK3}
Let $G$ be a graph on $n$ vertices such that each vertex is contained in at least $\binom t2$ triangles, where $t$ is a positive real number. How many edges must $G$ have?
\end{problem}

Our Theorem~\ref{thm:G} (asymptotically) answers Problem~\ref{pro:KK3} for $n/2\le t+1\le n$. 
The following proposition shows that for larger $n$, all but $O(t^3)$ vertices of an extremal graph are contained in isolated copies of $K_{t+1}$.

\begin{proposition}
When $n> (t+1)^2 (t+2)/4$, every extremal graph for Problem~\ref{pro:KK3} contains an isolated copy of $K_{t+1}$.
\end{proposition}
\begin{proof}
Let $G=(V, E)$ be an extremal graph with $|V|=n$. Since every vertex lies in at least $\binom t2$ triangles, it suffices to show that $G$ contains a vertex of degree $t$ and all of its neighbors also have degree $t$ (thus inducing an isolated copy of $K_{t+1}$).

Suppose $n= a(t+1) + b$, where $0< b\le t$. Let $G'$ be the disjoint union of $a-1$ copies of $K_{t+1}$ together with two copies of $K_{t+1}$ sharing $t+1-b$ vertices. 
Since $G$ is extremal, we have 
\[
2e(G)\le 2e(G')= t n + (t+1-b)b\le t n + (t+1)^2/4.
\]
Partition $V(G)$ into $A\cup B$ such that $A$ consists of all vertices of degree greater than $t$ and $B$ consists of all vertices of degree exactly $t$. Then 
\[
\sum_{v\in A} (d_G(v) - t) = \sum_{v\in V} (d_G(v) - t) = 2e(G) - tn \le (t+1)^2/4. 
\]
This implies that $|A|\le (t+1)^2/4$. Let $e(A, B)$ denote the number of edges (of $G$) between $A$ and $B$. It follows that 
\[
e(A, B)\le \sum_{v\in A} d(v) \le \frac14 (t+1)^2 + t|A|\le \frac14 (t+1)^3.
\]
Let $B'$ consists of the vertices of $B$ that are adjacent to some vertex of $A$. Then $|B'|\le e(A, B)\le (t+1)^3/4$. If $n> (t+1)^2 (t+2)/4$, then $n> |A| + |B'|$ and consequently, there exists a vertex of $B$ whose $t$ neighbors are all in $B$, as desired.
\end{proof}

The $t=2$ case of Problem~\ref{pro:KK3} assumes that every vertex in an $n$-vertex graph is contained in a triangle. Since $\d(G)\ge 2$, it follows that $e(G)\ge n$, which is best possible when $3$ divides $n$. Recently, Chakraborti and Loh \cite{MR4111669} determined the minimum number of edges an $n$-vertex graph in which every vertex is contained in a copy of $K_s$, for arbitrary $s\le n$. Their extremal graph is the union of copies of $K_s$, all but two of which are isolated.

\medskip
Finally, using careful case analysis, we can answer Problem~\ref{pro:KK3} \emph{exactly} when $t$ is very close to $n$. This falls into the $r+t> 5n/6$ case of Theorem~\ref{thm:G} but $G_2$ defined in Construction~\ref{cons}  is not necessarily extremal (unless both $r$ and $n/2$ are integers).
\begin{itemize}
\item When $n=t+2$, the (unique) extremal graph is $K_n^-$, the complete graph on $n$ vertices minus one edge.
\item When $n=t+3$ is even, the (unique) extremal graph is $K_n$ minus a perfect matching (provided $t>5$). When $n=t+3$ is odd, $K_n$ minus a matching of size $\frac{n-1}2$ is an extremal graph (provided $t>6$).
\item When $n=t+4$, the complement of any $K_3$-free 2-regular graph on $n$ vertices is an extremal graph. Note that $r= n/2 - 2$ in this case and thus $G_2$ is one of the extremal graphs when $n$ is even. 
\end{itemize}

\section*{Acknowledgment}
The authors thank Guantao Chen and Po-Shen Loh for valuable discussions.
We also thank two anonymous referees for their helpful comments.

\end{document}